\documentclass[leqno,12pt]{amsart}
 
\usepackage{amscd,amsfonts,amsmath,amssymb,amstext,amsthm}
\newtheorem{theorem} {Theorem} [section]
\newtheorem{lemma} [theorem]{Lemma}
\newtheorem{proposition}[theorem]{Proposition}
\newtheorem{corollary} [theorem]{Corollary}
\begin{document}
\title[Satake Diagrams and Real Structures]{Satake Diagrams\\
and\\
Real Structures on Spherical Varieties}
\author{Dmitri Akhiezer}
\address
{Dmitri Akhiezer\newline
Institute for Information Transmission Problems\newline
19 B.Karetny per.,127994 Moscow, Russia}
\noindent
\email{akhiezer@iitp.ru}
\begin{abstract}
With each antiholomorphic involution $\sigma $ of a connected complex
semisimple Lie group $G$ we associate an automorphism $\epsilon _\sigma$
of its Dynkin diagram.   
The definition of $\epsilon _\sigma $ is given in terms of the Satake diagram
of $\sigma $.
Let $H\subset G$ be a self-normalizing spherical subgroup. 
If $\epsilon _\sigma ={\rm id}$ then we prove
the uniqueness and existence of a $\sigma $-equivariant real structure
on $G/H$ and on the wonderful completion of $G/H$. 
\end{abstract}

\renewcommand{\subjclassname}
{\textup{2010} Mathematics Subject Classification}
\maketitle
\renewcommand{\thefootnote}{}
\footnotetext{Supported by SFB/TR 12 and SPP 1388 of the DFG}
\maketitle

\section{Introduction and statement of results}\label{intro}
In this paper, we consider real structures
on complex manifolds acted on by complex Lie groups.
A real structure on a complex manifold $X$
is an antiholomorphic involutive diffeomorphism $\mu : X \to X$.
Suppose a complex Lie group $G$ 
acts holomorphically on $X$ and let $\sigma : G \to G$ be an
involutive antiholomorphic
automorphism of $G$ as a real Lie group.
A real structure $\mu : X \to X$ is said to be $\sigma $-equivariant
if $\mu $ satisfies 
$\mu (g\cdot x) = \sigma (g) \cdot \mu (x)\ \ {\rm for\ all}\ g\in G, x\in X$.
We start with 
homogeneous manifolds of arbitrary complex Lie groups. In Section ~\ref{equiv} 
we prove that a $\sigma $-equivariant real structure on $X = G/H$
exists and is unique if $H$ is self-normalizing
and $\sigma (H)$ and $H$ are conjugate by an inner
automorphism of $G$. The conjugacy of $H$ and $\sigma (H)$
is also necessary for the existence of a $\sigma $-equivariant
real structure.

Assume $G$ is connected and semisimple and denote by
$\mathfrak g$ the Lie algebra of $G$.
In Section ~\ref{aut}, with any antiholomorphic involution $\sigma : G \to G$
we associate an automorphism class
$\epsilon _\sigma \in {\rm Aut (\mathfrak g)}/{\rm
Int (\mathfrak g)}$
acting on the Dynkin diagram in
the following way. 
We choose a Cartan subalgebra of the real form $\mathfrak g_0 \subset \mathfrak g$
and the root ordering as in the classical paper of I.Satake ~\cite{S}.
Let $\Pi_\bullet$ (resp. $\Pi_\circ$) be the set
of compact (resp. non-compact) simple roots, $\omega :  \Pi_\circ \to 
\Pi_\circ$
the involutory self-map associated with $\sigma $. Denote by 
$W_\bullet$ the subgroup of the Weyl group $W$ generated
by simple reflections $s_\alpha $, where $\alpha 
\in \Pi_\bullet$, and let
$w_\bullet$ be
the element of maximal length in $W_\bullet$. Then $\epsilon _\sigma (\alpha )=
-w_\bullet(\alpha )$ for $\alpha \in \Pi_\bullet $ and $\epsilon _\sigma
(\alpha ) = 
\omega (\alpha )$ for $\alpha \in \Pi _\circ$. 
On the Satake diagram,
$\epsilon _\sigma$ 
interchanges the white circles connected by two-pointed arrows
and permutes the black ones
as  
the outer automorphism of order 2 for compact algebras
$A_n (n\ge 2)$, $D_n$ ($n$ odd), $E_6$, and identically otherwise.

Let $B \subset G$ be a Borel subgroup.
Then $\sigma $ acts on the character group ${\mathcal X}(B)$ in
a natural way. Namely, $\sigma (B) = cBc^{-1}$ for some
$c \in G$ and, given $\lambda \in {\mathcal X}(B)$,  
the character  
$$B \ni b 
\mapsto \lambda ^\sigma (b) : = \overline {\lambda (c^{-1}\sigma (b)c)}$$
is in fact independent of $c$.
In Section ~\ref{action} we show that the arising 
action coincides with the one
given by $\epsilon _\sigma$. 

In Section ~\ref{sph} we consider equivariant real structures on 
homogeneous spherical spaces. It turns out that,
under some natural conditions on a spherical
subgroup $H \subset G$, the homogeneous space $G/H$ possesses
a $\sigma $-equivariant real structure. More precisely,
we have the following result.

\begin{theorem}\label{main1}
Assume 
$\epsilon _\sigma = {\rm id}$. Then
any spherical subgroup $H \subset G$ is conjugate to $\sigma (H)$
by an inner automorphism of $G$, i.e., $\sigma (H) = a H a^{-1}$
for some $a \in G$. The map
$$\mu _0 : G/H \to G/H,\ \ \mu_0(g\cdot H) := \sigma (g)\cdot a \cdot H,$$
is correctly defined, antiholomorphic and 
$\sigma $-equivariant.
Moreover, if the subgroup $H$ is self-normalizing
then: {\rm (i)} $\mu _0$ is involutive,
hence a $\sigma $-equivariant real structure on $G/H$; {\rm (ii)} such a
structure is unique.
\end{theorem}
In Section ~\ref{wonderful} we prove a similar theorem for 
wonderful varieties. Wonderful varieties
were introduced by D.Luna ~\cite {Lu}, and we recall 
their definition in Section ~\ref{wonderful}.
Wonderful varieties can be viewed as equivariant completions of spherical
varieties with certain properties. If such a completion
exists, it is unique. 
Furthermore, if $H$ is a self-normalizing spherical subgroup
of a semisimple group $G$ then,
by a result of F.Knop ~\cite{K3}, $G/H$ has a wonderful completion. 
 
\begin{theorem}\label{main2}
Let $H$ be a self-normalizing spherical subgroup of $G$
and let $X$ be the wonderful completion of $G/H$. If 
$\epsilon _\sigma = {\rm id}$
then there exists one and only one $\sigma $-equivariant real structure
$\mu:X \to X$.
\end{theorem}
                                                                  
\noindent
{\bf Remarks.}\ 1.\,Assume that $\sigma $ defines a split form of $G$. Then it is easily
seen that 
$\epsilon _\sigma  = {\rm id}$. In the split case
Theorems ~\ref{main1} and ~\ref{main2} 
are joint results with S.Cupit-Foutou ~\cite {A-CF}.
In this case, the 
$\sigma $-equivariant real structure on a wonderful variety $X$
is called {\sl canonical}. Assume in addition that $X$ is strict, i.e. all
stabilizers (and not just the principal one) are self-normalizing,
and equip $X$ with its canonical real structure.  
Then 
~\cite{A-CF} contains an estimate of the number of orbits of the connected
component $G^\sigma _0$ on the real part of $X$.

2.\,The results of this paper are related to 
the theory of (strongly) visible actions
introduced by T.Kobayashi, see e.g.
~\cite{Ko}. Namely, given a holomorphic action of a complex Lie group
on a complex manifold, its 
antiholomorphic diffeomorphisms are used in ~\cite{Ko}
to prove that the action of a real form is
(strongly) visible.

\section{Equivariant real structures}\label{equiv}
A real structure on a complex manifold
$X$ is an antiholomorphic involutive diffeomorphism
$\mu : X \to X$.
The set of fixed points $X^\mu $ of $\mu $
is called the real part of $X$ with respect to $\mu $.
If $X^\mu \ne \emptyset $ then $X^\mu $
is a closed real submanifold in $X$ and
$$\dim _{\mathbb R}(X^\mu ) = \dim _{\mathbb C}(X).$$ 
Suppose a complex Lie group $G$ 
acts holomorphically on $X$ and let $\sigma : G \to G$ be an
involutive antiholomorphic
automorphism of $G$ as a real Lie group.
The fixed point subgroup $G^\sigma $ is a real form of $G$.
A real structure $\mu : X \to X$ is said to be $\sigma $-equivariant
if 
$$\mu (g x) = \sigma (g) \cdot \mu (x)\ \ {\rm for\ all}\ g\in G, x\in X.$$
For such a structure the set $X^\mu $ is stable under $G^\sigma $.
We are interested in equivariant real structures on homogeneous manifolds
and on their equivariant embeddings. 

\begin{theorem}\label{conj}
 Let 
$G$ be a complex Lie group, 
let $\sigma : G \to G$ be an antiholomorphic involution, and
let $H \subset G$ be a 
closed complex Lie subgroup.
If there exists a $\sigma $-equivariant real structure on $X=G/H$
then $\sigma (H)$ and $H$ are conjugate by an inner automorphism of $G$.
Conversely, if $\sigma (H)$ and $H$ are conjugate and $H$ is
self-normalizing then a $\sigma $-equivariant real structure on $X$
exists and is unique.
\end{theorem}

\noindent
{\bf Proof.}\ Suppose first that $\mu : X \to X$ is a
$\sigma $-equivariant real structure. Let $x_0 = e\cdot H$
be the base point and let $\mu (x_0) = g_0\cdot H$. For $h \in H$
one has
$$\mu (x_0) = \mu (hx_0) = \sigma (h)\cdot \mu (x_0),$$
showing that $\sigma (H) \subset g_0Hg_0^{-1}$. 
To prove the opposite inclusion, observe that
$g\cdot\mu (x_0) = \mu (x_0)$ is equivalent to $\mu (\sigma (g)\cdot x_0)
= \mu (x_0)$. This implies 
$\sigma (g)\cdot x_0 = x_0$, so that $\sigma (g) \in H$ and 
$g\in \sigma (H)$, hence  $g_0Hg_0^{-1} \subset \sigma (H)$. 

To prove the converse,
assume that $H$ is self-normalizing and 
$$g_0Hg_0^{-1} = \sigma (H)$$
for some $g_0 \in G$.
Let $r_{g_0}$                     
be the right shift $g \mapsto gg_0$.       
We have a map $\mu : X \to X$,
correctly defined by
$\mu (g\cdot H) = \sigma (g)g_0\cdot H $. 
The commutative
diagram
\begin{displaymath}
\begin{array}{ccccc} 
G& {\buildrel \sigma \over \longrightarrow}&G&{\buildrel r_{g_0}\over
\longrightarrow}& G\\ 
\downarrow& & & &\downarrow\\ 
X =G/H& &{\buildrel\mu\over\longrightarrow}& & X=G/H \ ,\\
\end{array}
\end{displaymath}
where the vertical
arrows denote the canonical projection
$g\mapsto g\cdot H$, shows that the map $\mu $ is antiholomorphic.
It is also clear that $\mu $ is a $\sigma $-equivariant map, i.e.,
$\mu (gx) = \sigma (g)\cdot \mu(x)$ for all $g\in G$. 
Therefore $\mu ^2$ is an automorphism
of the homogeneous space $X$, i.e., $\mu ^2$
is a biholomorpic self-map of $X$ commuting with the $G$-action. 
Since $H$ is self-normalizing, we see that $\mu ^2 = {\rm id}$.
Thus $\mu $ is a $\sigma $-equivariant real structure on $X$.
If $\mu ^\prime $ is another such structure then $\mu \cdot \mu ^\prime$
is again an automorphism of $X = G/H$, 
so $\mu \cdot \mu ^\prime = {\rm id}$ and $
\mu^\prime = \mu $.  \hfill{$\square$}

\noindent
{\bf Example.} Let $B$ be a Borel subgroup of a semisimple complex 
Lie group $G$ and let
$X = G/B$ be the flag manifold of $G$. It follows from Theorem ~\ref{conj} that  
a $\sigma $-equivariant real structure $\mu : X \to X$
exists for any $\sigma: G \to G$. One has $X^\mu \ne \emptyset $
if and only if the minimal parabolic subgroup of $G^\sigma $ is solvable
or, equivalently, if the real form has no compact roots.  
\section{Automorphism $\epsilon_\sigma$}\label{aut}
Let $\mathfrak g$ be a complex semisimple
Lie algebra, $\mathfrak g_0$ a real form of $\mathfrak g$, and                 
$\sigma : \mathfrak g \to \mathfrak g$  
the corresponding antilinear involution. 
In this section we define the automorphism $\epsilon _\sigma $ of the Dynkin 
diagram of $\mathfrak g$, cf. ~\cite{A1, A2} and ~\cite{On}, $\S9$.
We start by recalling 
the notions of compact and non-compact roots, see e.g.
~\cite {OV}, Ch. 5.                                               

Let $\mathfrak g_0 = \mathfrak k + \mathfrak p$ be a Cartan decomposition.
The corresponding Cartan involution extends to $\mathfrak g = 
\mathfrak g_0 + i \cdot \mathfrak g_0$ as an automorphism
$\theta $ of the complex Lie algebra $\mathfrak g$. Clearly, 
$\theta ^2 = {\rm id}$ and $\sigma \cdot \theta = \theta \cdot \sigma$.
Pick a maximal abelian subspace $\mathfrak a \subset \mathfrak p$ and denote
by $\mathfrak m$ its centralizer in $\mathfrak k$. Let $\mathfrak h^+$ be
a maximal abelian subalgebra in $\mathfrak m$. Then   
$\mathfrak h = \mathfrak h^+ + \mathfrak a$ is a maximal abelian subalgebra in $\mathfrak g_0$
and any such subalgebra containing $\mathfrak a$ is of that form.
The Cartan subalgebra $\mathfrak t = \mathfrak h + i \cdot \mathfrak h \subset \mathfrak g$
is stable under $\theta $ and $\sigma $. On the dual space $\mathfrak t^*$,
we have the dual linear 
transformation                                               
$\theta ^{\rm T}$ and the dual antilinear transformation $\sigma ^{\rm T}$:
$$\theta ^{\rm T}(\gamma )(A) = \gamma (\theta A),\ 
\sigma ^{\rm T}(\gamma )(A) = \overline {\gamma (\sigma A)}\ \ \
(\ \gamma \in {\mathfrak t}^*,\ A\in {\mathfrak t}\ ).$$
Let $\Delta $ be the set of roots of $(\mathfrak g, \mathfrak t)$ 
and let $\Sigma $ 
be the set of roots of $\mathfrak g$ with respect to
$\mathfrak a \otimes {\mathbb C} = \mathfrak a + i\cdot \mathfrak a$.
Put $\mathfrak t_{\mathbb R} = i\cdot \mathfrak h^+ + \mathfrak a$. This is  
a maximal real subspace of $\mathfrak t$ on which all roots take
real values. Choose a basis $v_1,
\ldots, v_r, v_{r+1}, \ldots, v_l$ in $\mathfrak t_{\mathbb R}$
such that $v_1, \ldots, v_r$ form a basis of $\mathfrak a$ and introduce
the lexicographic ordering in the dual real vector spaces 
${\mathfrak t_{\mathbb R}}^*$
and ${\mathfrak a}^*$.               
Then $\Delta \subset {\mathfrak t_{\mathbb R}}^*$,
$\Sigma \subset {\mathfrak a}^*$, and
$\varrho (\Delta \cup \{0\}) = \Sigma \cup \{0\}$
under the restriction map $\varrho : {\mathfrak t_{\mathbb R}}^* \to {\mathfrak a}^*$.
Let $\Delta ^{\pm}, \Sigma ^{\pm}$ be the sets
of positive and negative roots
with respect to the chosen orderings, $\Pi \subset \Delta ^+,
\Theta \subset \Sigma^+$
the bases, $\Delta _{\bullet}       
=\{\alpha \in \Delta\ \vert \ \varrho (\alpha) = 0\},\
\Delta _\circ = \Delta \setminus \Delta _\bullet$. 
The roots from $\Delta _\bullet $ and $\Delta _\circ $
are called compact
and non-compact roots, respectively.
Let $\Delta _\bullet ^\pm = \Delta ^\pm
\cap \Delta _\bullet,
\Delta _\circ ^\pm =  \Delta ^\pm \cap \Delta _\circ,
\Pi_\bullet = \Pi \cap \Delta _\bullet$ and                      
$\Pi_\circ = \Pi \cap \Delta _\circ $.
One shows that $\Delta _\bullet $ is a root system with basis
$\Pi _\bullet $. Also, $\varrho (\Delta _\circ ^\pm) = \Sigma ^\pm,\ 
\theta ^{\rm T}(\Delta _\circ ^\pm) = \Delta _\circ ^\mp$ and
$\varrho (\Pi_\circ) = \Theta $.
Furthermore, one has an 
involutory self-map
$\omega : \Pi_\circ \to \Pi_\circ $,
defined by
$$\theta ^{\rm T}(\alpha) = -\omega(\alpha) - \sum_{\gamma \in
\Pi_\bullet }\  
c_{\alpha \gamma}\, \gamma,$$
where $c_{\alpha \gamma}$ are non-negative integers.
The Satake diagram is the Dynkin diagram on which the simple roots
from $\Pi _\bullet $ are denoted by black circles,
the simple roots from $\Pi_\circ$ by white circles, and
two white circles
are connected by a two-pointed arrow if and only if
they correspond to the roots $\alpha $ and $\omega (\alpha) \ne \alpha$.
 
Let $W$ be the Weyl group of $\mathfrak g$ with respect to $\mathfrak t$
considered as a linear group on
$\mathfrak t^*$.  
The subgroup of $W$ generated by the reflections $s_\alpha $ with $\alpha 
\in \Pi _\bullet $ is denoted by $W_\bullet $. The element of maximal
length in $W _\bullet $ with respect to these generators is denoted
by $w_\bullet$. Note that 
$-w_\bullet (\alpha ) \in \Pi_\bullet$ 
if $\alpha \in \Pi_\bullet $.
Let $\iota : \mathfrak g \to \mathfrak g$
be an inner automorphism such that $\iota (\mathfrak t) = {\mathfrak t}$,
acting as $w_\bullet $ on ${\mathfrak t}^*$. Since $w_\bullet ^2 = {\rm id}$,
we have
$$(\iota ^{\pm 1}\vert {\mathfrak t})^{\rm T} = w_\bullet.$$
Let $\eta $ be the Weyl involution of $\mathfrak g$ acting as
$-{\rm id}$ on $\mathfrak t$.   

\begin{proposition}\label{eps}
The self-map of $\Pi$, defined by
$$\epsilon_\sigma (\alpha ) = \begin{cases} 
-w_\bullet (\alpha ) & 
\mbox{if} \ \alpha \in \Pi_\bullet ,\\
\ \ \omega (\alpha ) &  \mbox{if} \ \alpha \in \Pi_\circ ,\\ 
\end{cases}$$
is an automorphism of the Dynkin diagram.
\end{proposition}

\noindent
{\bf Proof.}\ We must find an automorphism $\phi :\mathfrak g
\to \mathfrak g$ preserving $\mathfrak t$ and $\Pi $,
which acts on $\Pi $ as $\epsilon _\sigma $. 
Let $\phi = \eta 
\cdot \theta \cdot \iota $.
Then $\phi $ acts on $\Delta $ by
$$\alpha \mapsto - w_\bullet (\theta ^{\rm T}(\alpha )).$$
If $\alpha \in \Pi_\bullet $ then $\theta ^{\rm T}(\alpha ) = \alpha$,
and so $\phi $ sends $\alpha $ to 
$-w_\bullet (\alpha ) = \epsilon _\sigma (\alpha)$.
Now, if $\alpha  \in \Pi_\circ $ then
$$-w_\bullet (\theta ^{\rm T}(\alpha ))
=
w_\bullet (\omega (\alpha )) + \sum _{\gamma \in \Pi_\bullet}\ 
c_{\alpha \gamma}\, w_\bullet (\gamma)$$
by the definition of $\omega $. The simple reflections in
the decomposition of $w_\bullet $ correspond to the elements of 
$\Pi _\bullet$. Applying these
reflections to $\omega (\alpha )
\in \Pi_\circ $ one by one, we see that the right hand side
is the sum of $\omega (\alpha )$ and a linear combination of elements of 
$\Pi _\bullet $, whose coefficients must be nonnegative.
Therefore $-w_\bullet (\theta ^{\rm T}(\Pi )) \subset \Delta ^+$.
Since $-w_\bullet \cdot \theta^{\rm T}$ arises from $\phi $, this is
an automorphism of $\Delta $.
Thus $-w_\bullet (\theta ^{\rm T}(\Pi)) $ is a base of $\Delta $, hence
$-w_\bullet (\theta ^{\rm T}(\Pi)) = \Pi$. In particular,
$-w_\bullet (\theta ^{\rm T}(\alpha )) \in \Pi$,
and so we obtain $-w_\bullet (\theta ^{\rm T}(\alpha )) = \omega (\alpha ) = 
\epsilon _\sigma (\alpha )$.
\hfill $\square $

\begin{proposition}\label{eps1}
Extend $\epsilon _\sigma $ to a linear map of $\mathfrak t^*$
and denote the extension again by $\epsilon _\sigma$. Then
$w_\bullet $ and $\theta ^{\rm T}$ commute and
$$\epsilon _\sigma = -w_{\bullet }\theta ^{\rm T}
= -\theta ^{\rm T}w_\bullet.$$ 
\end{proposition}

\noindent
{\bf Proof.}
We already proved that $\epsilon _\sigma $ equals
$-w_{\bullet}\theta ^{\rm T}$ on $\Pi $, so it suffices to show
that $\epsilon _\sigma $ also equals $-\theta ^{\rm T}w_\bullet $
on $\Pi $. For $\alpha \in \Pi _\bullet $ we have $-w_\bullet (\alpha )
\in \Pi_\bullet $ and $\theta ^{\rm T}\alpha = \alpha $. Thus
$-w_\bullet \theta ^{\rm T}\alpha = - w_\bullet \alpha = 
-\theta^{\rm T}w_\bullet \alpha $.
For $\alpha \in \Pi_\circ $ we have
$$w_\bullet (\alpha )= \alpha +\sum _{\gamma \in \Pi_\bullet}
d_{\alpha \gamma}\ \gamma , \ \  d_{\alpha \gamma} \ge 0,$$
by the definition of $w_\bullet $. Applying $\theta ^{\rm T}$ we get
$$\theta ^{\rm T}w_\bullet (\alpha )= \theta ^{\rm T}(\alpha ) +\sum
_{\gamma \in \Pi_\bullet}\ d_{\alpha \gamma }\ \gamma = 
-\omega (\alpha ) - \sum _{\gamma \in \Pi_\bullet}\ (c_{\alpha \gamma}
- d_{\alpha \gamma})\ \gamma,$$
hence $-\theta ^{\rm T}w_\bullet (\alpha )\in \Delta ^+$. 
But 
$-\theta ^{\rm T}w_\bullet $ is an automorphism of $\Delta $.
Namely, define an automorphism $\phi ^\prime : \mathfrak g \to
\mathfrak g$ by $\phi ^\prime 
= \eta \cdot \iota \cdot \theta.$
Then $\phi ^\prime (\mathfrak t) = \mathfrak t$
and the dual to $\phi ^\prime \vert _{\mathfrak t}$ 
is $-\theta ^{\rm T}w_\bullet$.
Therefore
$-\theta ^{\rm T} w_\bullet (\Pi ) = \Pi $, so that $c_{\alpha \gamma}
= d _{\alpha \gamma}$ and
$-\theta ^{\rm T} w_\bullet (\alpha )= \omega (\alpha ) = 
\epsilon _\sigma (\alpha)$.
\hfill {$\square $}

\begin{proposition}\label{eps2}
Let ${\mathfrak b}^+$ be the positive Borel subalgebra defined by
the chosen ordering of roots, i.e.,
${\mathfrak b}^+ = 
{\mathfrak t}\ + \sum _{\alpha \in \Delta ^+}\ {\mathfrak g}_\alpha $.
Then                                                             
$$\sigma (\mathfrak b^+) = 
\iota (\mathfrak b^+) = \iota ^{-1}(\mathfrak b^+).$$
\end{proposition}

\noindent
{\bf Proof.} 
Observe that $\sigma \cdot \theta $ equals $-{\rm id}$ on
${\mathfrak t}_{\mathbb R}$.
Therefore $\sigma ^{\rm T}(\gamma ) = - \theta ^{\rm T}(\gamma )$
and $\sigma ^{\rm T} w_\bullet (\gamma ) = - \theta ^{\rm T} w_\bullet
(\gamma ) = \epsilon _\sigma (\gamma )$ for
$\gamma \in {\mathfrak t}_{\mathbb R}^*$. 
In particular,
for $\alpha \in \Delta ^+ $ we have
$\alpha ^\prime = \sigma ^{\rm T} w_\bullet (\alpha ) \in \Delta ^+$,
hence
$\sigma \cdot \iota ^{\pm 1}({\mathfrak g}_\alpha )= {\mathfrak g}_{\alpha ^\prime }
\subset \mathfrak b^+,$
and our assertion follows.
\hfill {$\square $}  

\noindent
{\bf Remark.}\ If $\mathfrak g$ is a complex simple Lie algebra considered
as a real one, then the Dynkin diagram
of its complexification is disconnected and has two isomorphic
connected components. Furthermore, $\Pi_\bullet = \emptyset$
and $\omega : \Pi_\circ \to \Pi_\circ$ maps each
component of the Satake diagram onto the other one. In particular,
$\epsilon _\sigma \ne {\rm id}$.   
If $\mathfrak g$ is simple and has no complex structure, then it is easy
to find the 
maps $\epsilon _\sigma $ for all Satake diagrams, see
~\cite{On}, Table 5. Let $l$ be the rank of $\mathfrak g$.
It turns out that $\epsilon _\sigma ={\rm id}$ for $\sigma $
defining ${\mathfrak sl}_{l+1}({\mathbb R}),\ 
{\mathfrak sl}_m({\mathbb H}),\ l=2m-1$, 
${\mathfrak so}_{p,q},\,
p+q = 2l, l\equiv p ({\rm mod}\,2)$, ${\mathfrak u}_l^*({\mathbb H}),\ l=2m$,
${\rm EI, EIV}$ or 
any real form of ${\rm B}_l, {\rm C}_l$,
${\rm E}_7, {\rm E}_8,\ {\rm F}_4$
and ${\rm G_2}$. For the remaining real forms $\epsilon _\sigma \ne {\rm id}$. 
\section{Action of $\sigma $ on ${\mathcal X}(B)$}\label{action}
In the rest of this paper $G$ is a connected complex semisimple Lie group.
Let 
$B\subset G$ be a Borel subgroup and $T \subset B$ a maximal torus. 
The Lie algebras are denoted by the corresponding German letters.
We want to apply the results of the previous section to
the automorphisms of $\mathfrak g$ which lift to $G$.
Suppose $\sigma $ is an antiholomorphic 
involutive automorphism of $G$ and denote again by $\sigma $  
the corresponding antilinear involution of $\mathfrak g$.
The automorphisms $\eta,\theta$ and $\iota $
lift to $G$ and the liftings are denoted by the same letters.
Recall that $\epsilon _\sigma $
is originally defined by its action on $\Pi$ as 
an automorphism class in ${\rm Aut}(\mathfrak g)/{\rm Int}(\mathfrak g)$.
The linear map induced by $\epsilon _\sigma $ on $\mathfrak t^*$
is denoted again by $\epsilon _\sigma$.
The automorphism 
$$\phi : \mathfrak g \to \mathfrak g, \ \  
\phi = \eta \cdot \theta \cdot \iota ,$$ 
leaves $\mathfrak t$ stable and acts
on $\mathfrak t^*$ as $\epsilon _\sigma $, see Propositions ~\ref{eps}
and ~\ref{eps1}. 
Since $\sigma $ and $\phi $ are globally defined,
we may consider their actions 
on the character groups of $T$ or $B$.

Namely, since $\sigma (B)$ is also a Borel subgroup, we have $\sigma (B) =
cBc^{-1}$ for some $c \in G$. 
The action of $\sigma $
on the character group ${\mathcal X}(B)$, given by
$$\lambda \mapsto \lambda ^\sigma,\ \ \ \lambda ^\sigma (b) = 
\overline {\lambda (c^{-1}\sigma (b) c)}\ \quad \ (b\in B), $$  
is correctly defined. For, if $d\in G$
is another element such that $\sigma (B) = dBd^{-1}$ then
$d^{-1}c \in B$, hence
$\lambda (d^{-1}\sigma (b) d) =
\lambda (d^{-1}c) \lambda (c^{-1}\sigma (b) c)\lambda (c^{-1}d) =
\lambda (c^{-1}\sigma(b)c)$.

Also, 
we have the right action of the automorphism group ${\rm Aut}(G)$
on ${\mathcal X}(B)$, defined in the same way. Namely,
for an automorphism $\varphi :G \to G$ we put 
$$\lambda ^\varphi (b) = \lambda (c^{-1} \varphi (b) c)\ \quad \ (b\in B),$$
where $c$ is chosen so that $\varphi (B) = cBc^{-1}$.

For two Borel subgroups $B_1, B_2$ the character groups
are canonically isomorphic. Moreover,
if $\lambda _1 \in {\mathcal X}(B_1)$
corresponds to $\lambda _2 \in {\mathcal X}(B_2)$ under the canonical
isomorphism then $\lambda _1^\sigma $ corresponds to $\lambda _2 ^\sigma $
and $\lambda _1 ^\varphi $ corresponds to $\lambda _2^\varphi $.

Clearly, $\lambda ^\varphi = \lambda$ for $\varphi \in {\rm Int}(G)$,
so we obtain the action
of ${\rm Aut}(G)/{\rm Int}(G)$ 
on ${\mathcal X}(B)$. In particular, we write 
${\epsilon _\sigma}(\lambda)$ instead of $\lambda ^\phi $.   

\begin{lemma}\label{lambda}
For any $\lambda \in {\mathcal X}(B)$ one has
$$\lambda^\sigma = {\epsilon_\sigma}(\lambda).$$ 
\end{lemma}

\noindent
{\bf Proof.}\ Choose $\mathfrak t $ and 
$\mathfrak b = \mathfrak b^+$ as in Section ~\ref{aut}. Then
$\sigma(B) = \iota(B)$ by Proposition ~\ref{eps2}.
Let $d\lambda$ be the differential of 
a character $\lambda $ at the neutral point of $T$.
Since
$\lambda ^\sigma (t) = \overline{\lambda (\iota ^{-1}\sigma (t))}$
for $t\in T$,
we have $d\lambda ^\sigma = \sigma^{\rm T}w_\bullet d\lambda $.
On the other hand,
${\epsilon _\sigma} (\lambda ) = \lambda ^\phi $, where 
$\phi = \eta \cdot \theta \cdot \iota$. 
In the course of the proof of Proposition ~\ref{eps} we have shown
that $\phi $ preserves $\mathfrak b^+$.  
Thus
$${\epsilon _\sigma}(\lambda )(t) = 
\lambda(\eta \theta \iota (t)) =
\lambda(\theta \iota (t))^{-1} \ \ 
(t\in T),$$ 
hence
$d{\epsilon _\sigma}(\lambda)  
=- w_\bullet \theta ^{\rm T} d\lambda = - \theta ^{\rm T}w_\bullet d\lambda $
by Proposition ~\ref{eps1}.
Since $\theta ^{\rm T} = - \sigma ^{\rm T}$
on $\mathfrak t^*_{\mathbb R}$, 
it follows that $d{\epsilon _\sigma}(\lambda) = d\lambda ^\sigma $.
\hfill {$\square $}

\noindent
{\bf Remark.}\ The automorphism class $\epsilon _\sigma $
has the following meaning for the representation theory, see ~\cite{A2}.
Let $V$ be an irreducible $G$-module with highest weight $\lambda $.
Denote by $\overline{V}$ the complex dual to the space of antilinear
functionals on $V$. Then $G$ acts on $\overline{V}$ in a natural way,
the action being antiholomorphic. This action combined with $\sigma $
is then holomorphic, the corresponding $G$-module is irreducible and
has highest weight
${\epsilon_\sigma}(\lambda)$.

\section{Spherical homogeneous spaces}\label{sph}
Let
$X = G/H$ be a spherical homogeneous space.  
We fix a Borel subgroup $B \subset G$                            
and recall the definitions of Luna-Vust invariants of $X$, see ~\cite{LV, K1,T}. 

For $\chi \in {\mathcal X}(B)$ let $^{(B)}{\mathbb C}(X)_\chi \subset {\mathbb C}(X)$
be the subspace 
of rational $B$-eigenfunctions of weight $\chi $, i.e.,
$$^{(B)}{\mathbb C}(X)_\chi = \{f\in {\mathbb C}(X)\ \vert \ f(b^{-1}x)
= \chi (b)f(x)\ \ (b\in B,x\in X)\}. $$ 
Since $X$ has an open $B$-orbit, this subspace is 
either trivial or one-dimensional. In the latter case we choose
a non-zero function $f_\chi \in \ ^{(B)}{\mathbb C}(X)_\chi $.  
The weight lattice $\Lambda (X)$ is the set of $B$-weights in ${\mathbb C}(X)$,
i.e.,
$${\Lambda }(X) 
= \{\chi \in {\mathcal X}(B)\ \vert \ ^{(B)}{\mathbb C}(X) _\chi \ne \{0\}\}.$$
Let ${\mathcal V}(X)$ denote the set of $G$-invariant discrete 
${\mathbb Q}$-valued
valuations of ${\mathbb C}(X)$. The
mapping
$${\mathcal V}(X) \to {\rm Hom}(\Lambda (X), {\mathbb Q}),\ \ v
\mapsto \{\chi \to  v(f_\chi)\}$$
is injective, see ~\cite{LV,K3},
 and so we regard ${\mathcal V}(X) $ as a subset of
$${\rm Hom}(\Lambda (X), {\mathbb Q}).$$ 
It is known that ${\mathcal V}(X)$
is a simplicial cone, see ~\cite{B,K2}. 

The set of all $B$-stable prime divisors in $X$ is denoted by ${\mathcal D}(X)$.
This is a finite set.
To any $D\in {\mathcal D}(X)$ 
we assign  $\omega _D \in {\rm Hom}(\Lambda (X), {\mathbb Q})$. Namely,
$\omega _D (\chi) = {\rm ord}_D f_\chi$, the order of $f_\chi $ along $D$.
We also write $G_D$ for the stabilizer of $D$.
The Luna-Vust invariants of $X$ are given by the
triple $\Lambda (X), {\mathcal V}(X), {\mathcal D}(X)$.
The homogeneous space $X$ is completely determined by these
combinatorial invariants. More precisely, one has the following theorem
of I.Losev ~\cite{Lo}.

\medskip

\begin{theorem}~\label{Losev} 
Let $X_1= G/H_1, X_2= G/H_2$
be two spherical homogeneous spaces. Assume that 
${\Lambda}(X_1) = {\Lambda }(X_2),\ 
{\mathcal V}(X_1) = {\mathcal V}(X_2)$. Assume further there is a bijection
$j: {\mathcal D}(X_1) \to {\mathcal D}(X_2)$, such that
$\omega_D = \omega_{j(D)}, G_D = G_{j(D)}$. Then $H_1$
and $H_2$ are conjugate by an inner automorphism of $G$. 
\end{theorem}

\noindent
We now return to equivariant real structures. Let $\sigma $ be an 
antiholomorphic involution of a semisimple complex algebraic group.
Given a spherical subgroup $H \subset G$, observe that $\sigma (H)$
is also a spherical subgroup of $G$. Put $X_1 = G/H, X_2 = G/\sigma (H)$,
and denote again by $\sigma $ the antiholomorphic map 
$$X_1 \to X_2,\ g\cdot H \mapsto \sigma(g)\cdot \sigma (H).$$
Since 
the conjugate coordinate functions of $\sigma : G \to G$ are regular,
we have
$\sigma ^*\cdot {\mathbb C}(X_2) = \overline {{\mathbb C}(X_1)}$. 
Choose and fix $c\in G$
in such a way that $\sigma (B) = cBc^{-1}$. 

\begin{proposition}\label{lat}
$\epsilon _\sigma (\Lambda (X_1)) = 
\Lambda (X_2)$.
\end{proposition}

\noindent
{\bf Proof.}
For $f \in {\mathbb C}(X_2)$ define a rational function on $X_1$ by
$$f^\prime (x) = \overline{f(\sigma (cx))}.$$
Note that for $b \in B$ one has
$b^\prime := \sigma (cbc^{-1})\in B$.
Furthermore, since $b_0:=\sigma (c)c\in B$, we have
$$\chi ^\sigma (b) = 
\overline {\chi (c^{-1}\sigma (b) c)} =
\overline {\chi (b_0^{-1}\sigma (c)\sigma (b)\sigma (c)^{-1}b_0)} =
\overline {\chi (b^\prime)}.$$
Now take $f = f_\chi $. Then we obtain
$$f^\prime(b^{-1}x) =
\overline{f(\sigma(cb^{-1}x))}=
\overline{f(\sigma({b^\prime}^{-1})\sigma (cx))}
=
 \overline{\chi (b^\prime )}f^\prime (x),$$
showing that $f^\prime $ is a $B$-eigenfunction 
of weight $\chi ^\sigma $ on $X_1$.
Since the transform $f\mapsto f^\prime $ is invertible and $\chi ^\sigma =
\epsilon _\sigma (\chi )$ by 
Lemma ~\ref{lambda}, it follows that $\Lambda (X_2) = 
\epsilon _\sigma (\Lambda (X_1))$.
\hfill $\square $ 

\begin{proposition}\label{norm}
Extend $\epsilon _\sigma $ by duality to
${\rm Hom}({\mathcal X}(B),{\mathbb Q})$. Then $
\epsilon _\sigma ({\mathcal V}(X_1)) = {\mathcal V}(X_2)$. 
\end{proposition}

\noindent 
{\bf Proof.} 
The map
$${\mathbb C}(X_2)\ni f \mapsto \overline{f\circ \sigma} \in {\mathbb C}(X_1)$$
is a field isomorphism which is $\sigma $-equivariant in the obvious sense,
namely,
$$\overline {(g\cdot f)\circ \sigma} = \sigma(g)\cdot 
(\overline {f\circ \sigma})\
\ (g\in G).$$
Therefore, for $v \in {\mathcal V}(X_1)$ the valuation of ${\mathbb C}(X_2)$
defined by $v^\prime (f) = v(\overline {f\circ \sigma })$
is also $G$-invariant, i.e., $v^\prime \in {\mathcal V}(X_2)$.
Furthermore, since the function $f^\prime $,
defined in Proposition ~\ref{lat}, is in the $G$-orbit of 
$\overline {f\circ \sigma}$, we have
$v^\prime (f) = v(f^\prime )$. Now take $f=f_\chi $. Then $f^\prime $
is a $B$-eigenfunction with weight $\epsilon _\sigma (\chi )$.
Therefore
$\epsilon _\sigma (v) = v^\prime $.  
\hfill $\square $

For a $B$-invariant divisor $D$ on $X_1$ its image $\sigma (D)$
is a $\sigma (B)$-invariant divisor on $X_2$. Obviously, the map
$$j: {\mathcal D}(X_1) \to {\mathcal D}(X_2),\                             
j(D) := \sigma (c\cdot D),$$
is 
a bijection.

\begin{proposition}\label{div}
For any $D\in {\mathcal D}(X_1)$ one has
$\omega _{j(D)}=\epsilon _\sigma (\omega _D)$.
The stabilizers of $D$ and $j(D)$
are parabolic subgroups containing $B$ and satisfying
$$\sigma (G_{j(D)}) = cG_Dc^{-1}.$$
Their root systems are obtained from each other by $\epsilon _\sigma$.
\end{proposition} 

\noindent
{\bf Proof.}\ 
Let $f \in {\mathbb C}(X_2)$ and let $f^\prime \in {\mathbb C}(X_1)$ be
the function defined in 
Proposition ~\ref{lat}. Then 
$${\rm ord}_{j(D)}\, f = {\rm ord}_D\, f^\prime .$$
Applying this to $f = f_\chi $ we obtain $\omega _{j(D)} = \epsilon _\sigma 
(\omega _D)$.
The definition of $j$ implies readily
that $\sigma (G_{j(D)}) = cG_Dc^{-1}$, and 
the last assertion follows from 
Lemma ~\ref{lambda}
\hfill $\square $

\noindent
Combining Propositions ~\ref{lat}, ~\ref{norm}, and ~\ref{div}, 
we get the following
corollary.

\begin{corollary} 
If 
$\epsilon _\sigma $ leaves stable $\Lambda (X_1), {\mathcal V}(X_1)$
and, for any $D \in {\mathcal D}(X_1)$, one has $\epsilon _\sigma (\omega _D) =
\omega _D$ and $\sigma (G_D) = cG_Dc^{-1}$ then $H$ and $\sigma (H)$
are conjugate by an inner automorphism, i.e., $\sigma (H) = a H a^{-1}$,
where $a \in G$.
The map
$g\cdot H \mapsto \sigma (g)a \cdot H$
is correctly defined, antiholomorphic and 
$\sigma $-equivariant.
Moreover, if the subgroup $H$ is self-normalizing
then this map is 
a $\sigma $-equivariant real structure on $X_1$ and such a
structure is unique.
\end{corollary}

\noindent
{\bf Proof.} The conjugacy of $H$ and $\sigma (H)$
results from Theorem ~\ref{Losev}. The remaining
assertions follow from Theorem ~\ref{conj}.
\hfill $\square $ 

\medskip
\noindent
{\bf Proof of Theorem~\ref{main1}}. It suffices to apply
the above corollary in the case $\epsilon _\sigma ={\rm id}$.
\hfill $\square $


\begin{proposition}\label{inv}
If $\epsilon _\sigma =
{\rm id}$ and $\mu _0$ is defined as in Theorem ~\ref{main1}, then
any $v \in {\mathcal V}(G/H)$ is $\mu _0 $-invariant, i.e.,
for a non-zero rational function $f \in {\mathbb C}(G/H)$ one has
$v(\overline {f\circ \mu_0}) = v(f)$.
\end{proposition}

\noindent
{\bf Proof.}\ Consider $f$ as a right $H$-invariant function on $G$ and
put $f^a(g) = f(ga)\ (g\in G)$. Then $f^a$ is right $a H a^{-1}$-invariant.
Since $\sigma(H) = a H a^{-1}$, we can view $f^a$ as a rational function on
$X_2 = G/\sigma (H)$. Recall that we have the map $\sigma : X_1 \to X_2$.
The definition of $\mu _0 : X_1 \to X_1$ implies
$f\circ \mu_0 = f^a \circ \sigma$.
It suffices to prove the equality 
$v(\overline {f\circ \mu_0}) = v(f)$
on $B$-eigenfunctions.
Now, if $f$ is such a function then
$f^a$ is also a $B$-eigenfunction with the same weight. In the proof
of Proposition ~\ref{norm}, for a given $v
\in {\mathcal V}(X_1)$ we defined $v^\prime \in {\mathcal V}(X_2)$ and proved
that $\epsilon _\sigma (v) = v^\prime $. In our setting $v = v^\prime$,
and so we obtain $v(f) = v^\prime (f^a) = v(\overline {f^a\circ \sigma})=
v(\overline {f\circ \mu_0})$.
\hfill $\square $ 

\noindent
{\bf Example.}
Up to an automorphism of $X = {\mathbb C}{\mathbb P}^d$, there are two 
real structures $\mu_1, \mu_2: X \to X$ for $d$ odd and one real structure
$\mu_1: X \to X$ for $d$ even. 
In homogeneous coordinates 
$$\mu _1 (z_0:z_1:\ldots :z_d) =
(\overline {z_0}:\overline {z_1}: \ldots :\overline {z_d}) $$
and
$$\mu _2 (z_0:z_1:\ldots :z_d) = (- \overline {z_1}: \overline {z_0}:
\ldots :-\overline{z_d}:\overline{z_{d-1}}),\ d=2l-1.$$
One has $X^{\mu _1} = {\mathbb R}{\mathbb P}^d$
and $X^{\mu _2} = \emptyset$. 
Let $s_l$   
be the block $(2l\times 2l)$-matrix with $l$ diagonal blocks
\begin{displaymath}
\left(\begin{array}{cc} 0 & -1 \\ 1 & 0 \\ \end{array} \right).
\end{displaymath}
For 
$g\in G = {\rm SL}(d+1, {\mathbb C})$
put 
$$\sigma _1 (g)= \overline {g} \ \ {\rm and}\ \  
\sigma _2(g) = -s_l\overline {g}s_l\ \ {\rm if}\ \  d+1 = 2l.$$
Then
$G^{\sigma _1} = {\rm SL}(d+1, {\mathbb R})$ (the split real form) and
$G^{\sigma _2} = {\rm SL}(l, {\mathbb H})$, where $d+1 = 2l$.
One checks easily that $\mu _1$ is $\sigma _1$-equivariant
and $\mu _2$ is $\sigma _2$-equivariant. Note that a real structure
can be $\sigma $-equivariant only for one involution $\sigma $.
Therefore $X$ has
no $\sigma $-equivariant real structure if 
$\sigma$ defines a pseudo-unitary group ${\rm SU}(p,q),\ p+q = d+1$.
\section{Wonderful embeddings}\label{wonderful}
A complete non-singular algebraic $G$-variety $X$
of a semisimple group $G$ is called {\it wonderful}
if $X$ admits an open $G$-orbit whose
complement is a finite union of smooth prime divisors $X_1,\ldots,X_r$ with
normal crossings and the closures of $G$-orbits on $X$ are precisely the
partial intersections of these divisors. The notion
of a wonderful variety was introduced by D.Luna ~\cite{Lu},
who also proved that wonderful varieties are spherical. 
The total number of $G$-orbits on $X$ is $2^r$. The number
$r$ coincides with the rank of $X$ as a spherical variety. Moreover,
if a spherical homogeneous space $G/H$ has a wonderful embedding then such
an embedding is
unique up to a $G$-isomorphism.

\begin{theorem}\label{wond1} 
Let $G$ be a complex semisimple algebraic group,
$H\subset G$ a spherical subgroup, and $\sigma : G \to G $ an antiholomorphic
involution. Assume that $G/H$ admits a wonderful embedding 
$G/H \hookrightarrow X$. If there exists a $\sigma $-equivariant
real structure on $G/H$ then it 
extends to a $\sigma $-equivariant real structure
on $X$.
\end{theorem}

\noindent
{\bf Proof.}
This follows from the uniqueness of wonderful embedding.
Namely, let $\varepsilon: G/H \to X$ be the given wonderful embedding.
Take the complex conjugate $\overline {X}$ of $X$ and let
$\overline{\varepsilon}: G/H \to \overline{X}$ be the corresponding
antiholomorphic map.
We identify $\overline {X}$ with $X$ as topological spaces and endow
$\overline{X}$ with the action $(g,x) \mapsto \sigma(g)\cdot x$, which
is regular. Now, take a $\sigma $-equivariant
real structure $\mu _0$ on $G/H$ and consider the map $\overline
{\varepsilon} \circ \mu _0:
G/H \to \overline {X}$. This is again a wonderful embedding
of $G/H$. Since two wonderful embeddings are $G$-isomorphic,
there is a $G$-isomorphism $\nu : X \to \overline{X}$ such that
$\nu \circ \varepsilon = \overline{\varepsilon}\circ \mu _0$. 
The map $\nu $ defines
a required $\sigma $-equivariant real structure on $X$.
\hfill{$\square $}

\medskip
\noindent
{\bf Proof of Theorem ~\ref{main2}.}
Let $G/H \hookrightarrow X$ be the wonderful completion.
The existence and uniqueness
of a $\sigma$-equivariant real structure $\mu _0$ on $G/H $
follows from Theorem ~\ref{main1}. By Theorem ~\ref{wond1} 
this real structure extends to $X$,
the extension being obviously unique. \hfill {$\square $} 

\medskip
\noindent
As an application of our previous results we have the following
property of the $\sigma $-equivariant real structure $\mu $.

\begin{theorem} 
We keep the notations and assumptions of Theorem ~\ref{main2}.
Then all
$G$-orbits on $X$ are $\mu $-stable.
\end{theorem}
        
\noindent
{\bf Proof.}
It suffices to show that all divisors $X_i$ are $\mu $-stable.
Each $X_i$ defines a $G$-invariant valuation
of the field ${\mathbb C}(X) = {\mathbb C}(G/H)$. By Proposition ~\ref{inv}
such a valuation is $\mu $-invariant. Since the divisor is uniquely
determined by its valuation, it follows that $X_i$ are $\mu $-stable.
\hfill $\square $

\begin{corollary} Keeping the above notations and assumptions,
suppose that $\mu $ has a fixed
point in the closed $G$-orbit $X_1 \cap \ldots \cap X_r \subset X$.
Then $\mu $ has a fixed point in any $G$-orbit. In
particular, the number of $G^\sigma $-orbits in $X^\mu $ is
greater than or equal to $2^r$.
\end{corollary}

\noindent
{\bf Proof.} The closure of any $G$-orbit in $X$ is of the form 
$Y = X_{i_1}\cap \ldots \cap X_{i_k}$. 
We know that $Y$ is $\mu $-stable and has
a non-trivial intersection with $X^\mu $. Since the real
dimension of $X^\mu \cap Y$ equals the complex dimension of
$Y$, the set $X^\mu $ must intersect the open $G$-orbit in $Y$.
\hfill $\square $

\medskip
\noindent
The condition $\epsilon_\sigma = {\rm id}$ is essential.

\noindent
{\bf Example.} The adjoint representation of ${\rm SL}(2,{\mathbb C)}$
gives rise to a two-orbit action on the projective plane. The closed
orbit is the quadric $Q \subset {\mathbb C}{\mathbb P}^2$
arising from the nilpotent cone in the Lie algebra $\mathfrak{s}\mathfrak{l}(2,
{\mathbb C})$.
Let $G = {\rm SL}(2,{\mathbb C})\times {\rm SL}(2, {\mathbb C})$ and
$\sigma (g_1,g_2) =(\bar g_2, \bar g_1)$, 
where $g_1,g_2 \in {\rm SL}(2,\mathbb C)$. Note that $G^\sigma 
= {\rm SL}(2,\mathbb C)$ considered as a real group and $\epsilon _\sigma
\ne {\rm id}$. Let
$X = {\mathbb C}{\mathbb P}^2 \times                       
{\mathbb C}{\mathbb P}^2$ with each simple factor of $G$ acting
on the corresponding factor of $X$ in the way described above.
Then $X$ is a wonderful variety of rank 2. The divisors $X_1, X_2$
from the definition of a wonderful variety are ${\mathbb C}{\mathbb P}^2
\times Q$ and $Q \times {\mathbb C}{\mathbb P}^2$. The $\sigma $-equivariant
real structure $\mu $ on $X$ is given by $\mu (z_1,z_2) = (\bar z_2, \bar z_1),
\ \ 
z_1, z_2 \in {\mathbb C}{\mathbb P}^2$. The $G$-stable hypersurfaces $X_1,X_2$
are interchanged by $\mu $ and not $\mu$-stable.

\end {document}